# Hypothesis Testing with the General Source [†]

Te Sun HAN[‡]

April 26, 2000

[‡]Te Sun Han is with the Graduate School of Information Systems, University of Electro-Communications, Chofugaoka 1-5-1, Chofu, Tokyo 182-8585, Japan. E-mail: han@is.uec.ac.jp

**Abstract:** The asymptotically optimal hypothesis testing problem with the *general* sources as the null and alternative hypotheses is studied under exponential-type error constraints on the first kind of error probability. Our fundamental philosophy in doing so is first to convert all of the hypothesis testing problems completely to the pertinent computation problems in the *large deviation-probability theory*. It turns out that this kind of methodologically new approach enables us to establish quite compact general formulas of the optimal exponents of the second kind of error and correct testing probabbilities for the *general* sources including all nonstationary and/or non-ergodic sources with *arbitrary* abstract alphabet (countable or uncountable). Such general formulas are presented from the *information-spectrum* point of view.

# 1 Introduction

The hypothesis testing problem is very important not only from the theoretical viewpoint but also from the engineering point of view. This fundamental research subject in the hypothesis testing problem seems to have started earlier in the 1930's with the nonasymptotic study on that for one shot sources with *real* alphabet (e.g., see Neyman and Pearson [11]) and subsequently has been generalized into various kinds of directions including that of the asymptotic approach to a diversity of source processes.

In the present paper we consider a wide class of *general* sources (generalized processes) as null and alternative hypotheses. Let us first define the *general* source as an infinite sequence $\mathbf{X} = \{X^n = (X_1^{(n)}, \cdots, X_n^{(n)})\}_{n=1}^{\infty}$ of $n$-dimensional random variables $X^n$ where each component random variable $X_i^{(n)}$ ($1 \le i \le n$) takes values in an arbitrary *abstract* set $\mathcal{X}$ that we call the *source alphabet* (cf. Han [20]). It should be noted here that each component of $X^n$ may change depending on block length $n$. This implies that the sequence $\mathbf{X}$ is quite general in the sense that it may not satisfy even the consistency condition as usual processes, where the consistency condition means that for any integers $m, n$ such that $m < n$ it holds that $X_i^{(m)} \equiv X_i^{(n)}$ for all $i = 1, 2, \cdots, m$. The class of sources thus defined covers a very wide range of sources including all nonstationary and/or nonergodic sources. The introduction of such a class of general sources is crucial in the whole argument in the sequel. Thus, given two arbitrary *general* sources $\mathbf{X} = \{X^n\}_{n=1}^{\infty}$ and $\overline{\mathbf{X}} = \{\overline{X}^n\}_{n=1}^{\infty}$ taking values in the same source alphabet $\{\mathcal{X}^n\}_{n=1}^{\infty}$, we may define the general hypothesis testing problem with $\mathbf{X} = \{X^n\}_{n=1}^{\infty}$ as the *null hypothesis* and $\overline{\mathbf{X}} = \{\overline{X}^n\}_{n=1}^{\infty}$ as the *alternative hypothesis*.

**Remark 1.1** A more reasonable definition of the general source is the following. Let $\{\mathcal{Z}_n\}_{n=1}^{\infty}$ be any sequence of *arbitrary* source alphabets $\mathcal{Z}_n$ and let $Z_n$ be any random variable taking values in $\mathcal{Z}_n$ ($n = 1, 2, \cdots$). Then, the sequence $\mathbf{Z} = \{Z_n\}_{n=1}^{\infty}$ of random variables $Z_n$ is called a *general source*. The above definition is a special case of this general source with $\mathcal{Z}_n = \mathcal{X}^n$ ($n = 1, 2, \cdots$). The key results in this paper (Theorem 2.1 and Theorem 4.1) continue to be valid as well also in this more general setting with $\{\mathcal{X}^n\}_{n=1}^{\infty}$ (source alphabet), $\mathbf{X} = \{X^n\}_{n=1}^{\infty}$ (null hypothesis), $\overline{\mathbf{X}} = \{\overline{X}^n\}_{n=1}^{\infty}$ (alternative hypothesis) replaced by $\{\mathcal{Z}_n\}_{n=1}^{\infty}$ (source alphabet), $\mathbf{Z} = \{Z_n\}_{n=1}^{\infty}$ (null hypothesis), $\overline{\mathbf{Z}} = \{\overline{Z}_n\}_{n=1}^{\infty}$ (alternative hypothesis), respectively, where both of $Z_n$ and $\overline{Z}_n$ take values in $\mathcal{Z}_n$ ($n = 1, 2, \cdots$). □

In the present paper, with this kind of general hypothesis testings we investigate the optimal exponent problem for the probability of testing error as well as the optimal exponent problem for the probability of correct testing. Formally, let $\mathcal{A}_n$ be any subset of $\mathcal{X}^n$ ($n = 1, 2, \cdots$) that we call the *acceptance region* of the hypothesis testing, and define
$$\mu_n \equiv \Pr\{X^n \notin \mathcal{A}_n\}, \quad \lambda_n \equiv \Pr\{\overline{X}^n \in \mathcal{A}_n\}, \tag{1.1}$$
where $\mu_n, \lambda_n$ are called the *first kind of error probability* and the *second kind of error probability*, respectively.

One of the basic problems in the hypothesis testing is to determine the supremum $B_e(r|\mathbf{X}||\overline{\mathbf{X}})$ of achievable exponents for the second kind of *error* probability $\lambda_n$ under asymptotic constraints of the form $\mu_n \leq e^{-nr}$ on the first kind of error probability ($r > 0$ is a prescribed arbitrary constant) which means that the first kind of error probability is required at most to decay exponentially fast with the exponent $r$. Another basic problem in the hypothesis testing is to determine the infimum $B_e^*(r|\mathbf{X}||\overline{\mathbf{X}})$ of achievable exponents for the second kind of *correct* probability $1 - \lambda_n$ under asymptotic constraints of the same form as above $\mu_n \leq e^{-nr}$ on the first kind of error probability ($r > 0$ is again a prescribed arbitrary constant).

In the following sections we focus on these two basic problems for the general hypothesis testings. We establish a general formula (Theorem 2.1) for $B_e(r|\mathbf{X}||\overline{\mathbf{X}})$ in Section 2 along with several typical examples in Section 3, whereas we establish a general formula (Theorem 4.1) for $B_e^*(r|\mathbf{X}||\overline{\mathbf{X}})$ in Section 4 along with several typical examples in Section 5. In order to drive the general formula for $B_e(r|\mathbf{X}||\overline{\mathbf{X}})$ as well as that for $B_e^*(r|\mathbf{X}||\overline{\mathbf{X}})$ in a surprisingly unifying way, we shall take an *information-spectrum approach* that had been effectively invoked already in Han and Verdú [1], Verdú and Han [5], Han [17, 19, 20], where the substantially novel technique of *information spectrum slicing*, as exploited in Han [17, 18], plays the key role. Our fundamental philosophy here is first to convert all of the hypothesis testing problems completely to the pertinent computation problems in the *large deviation-probability theory*. We can then expel all the *acceptance-region* arguments from the original hypothesis testing problems; thereby, all of what we should do boils down solely to how to compute the relevant large deviation probabilities (or, in many standard cases, the relevant rate functions). It turns out that this kind of methodologically new approach enables us to establish quite compact general formulas of the exponent functions $B_e(r|\mathbf{X}||\overline{\mathbf{X}})$, $B_e^*(r|\mathbf{X}||\overline{\mathbf{X}})$ for general sources including all nonstationary and/or nonergodic sources with *abstract* alphabet. Such general formulas are presented in this paper.

Finally, in Section 6 we pleasingly observe that all the arguments developed in Sections 2∼5 continue to be valid even if we replace the general alternative hypothesis $\overline{\mathbf{X}} = \{\overline{X}^n\}_{n=1}^{\infty}$ by any sequence $\overline{\mathbf{X}} = \{G_n\}_{n=1}^{\infty}$ of nonnegative measures (for example, *counting* measures; *not necessarily* probability measures), and as a consequence in Section 7 it is revealed that there exists an intrinsic one-to-one operational correspondence between the problem of so generalized hypothesis testings and the problem of general fixed-length source codings. As an illustrative case, it is shown in the case of *countably infinite* source alphabet $\mathcal{X}$ that the general formula of Han [20] for the infimum $R_e(r|\mathbf{X})$ of achievable coding rates under asymptotic constraints of the form $\varepsilon_n \leq e^{-nr}$ ($r > 0$) on the error probability $\varepsilon_n$ with fixed-length source coding immediately follows from the general formula (Theorem 2.1) for $B_e(r|\mathbf{X}||\overline{\mathbf{X}})$ (with the sequence $\overline{\mathbf{X}} = \{C_n\}_{n=1}^{\infty}$ of counting measures) as driven in Section 2. It thus turns out that the general fixed-length source coding problem is just a special case of the so generalized hypothesis testing problem.

## 2  Hypothesis Testing and Large Deviation: Probability of Testing Error

In this section we investigate the problem of determining the supremum $B_e(r|\mathbf{X}||\overline{\mathbf{X}})$ of achievable exponents for the second kind of *error* probability $\lambda_n$ under asymptotic constraints of the form $\mu_n \leq e^{-nr}$ on the first kind of error probability $\mu_n$ ($r > 0$ is a prescribed arbitrary constant). Let us first give the formal definititons, where $\mathbf{X} = \{X^n\}_{n=1}^{\infty}$, $\overline{\mathbf{X}} = \{\overline{X}^n\}_{n=1}^{\infty}$ indicate the null hypothesis and the alternative hypothesis, respectively.

**Definition 2.1** A rate $E$ is called *r- achievable* if there exists an acceptance region $\mathcal{A}_n$ such that

$$\liminf_{n \to \infty} \frac{1}{n} \log \frac{1}{\mu_n} \geq r \quad \text{and} \quad \liminf_{n \to \infty} \frac{1}{n} \log \frac{1}{\lambda_n} \geq E.$$

**Definition 2.2 (The supremum of $r$-achievable error exponents)**

$$B_e(r|\mathbf{X}||\overline{\mathbf{X}}) = \sup \{E \mid E \text{ is } r\text{-achievable}\}.$$

The purpose of this section is to determine $B_e(r|\mathbf{X}||\overline{\mathbf{X}})$ as a function of $r$. To this end, we consider the random variable $\frac{1}{n} \log \frac{P_{X^n}(X^n)}{P_{\overline{X}^n}(X^n)}$ that we call the *divergence-density rate*,[*] and define the key function $\eta(R)$ by

$$\eta(R) = \liminf_{n \to \infty} \frac{1}{n} \log \frac{1}{\Pr\left\{\frac{1}{n} \log \frac{P_{X^n}(X^n)}{P_{\overline{X}^n}(X^n)} \leq R\right\}}, \tag{2.1}$$

where in the sequel we use the convention that $P_Z(\cdot)$ denotes the probability distribution of a random variable $Z$. It is obvious that this function $\eta(R)$ is monotone decreasing in $R$ but not necessarily continuous. Next, define the *spectral inf-divergence rate* $\underline{D}(\mathbf{X}||\overline{\mathbf{X}})$ of the random variable $\frac{1}{n} \log \frac{P_{X^n}(X^n)}{P_{\overline{X}^n}(X^n)}$ as

**Definition 2.3**
$$\underline{D}(\mathbf{X}||\overline{\mathbf{X}}) = \text{p-}\liminf_{n \to \infty} \frac{1}{n} \log \frac{P_{X^n}(X^n)}{P_{\overline{X}^n}(X^n)}.[\dagger]$$

**Lemma 2.1** If $R > \underline{D}(\mathbf{X}||\overline{\mathbf{X}})$, then $\eta(R) = 0$.

---

[*]In the case where the source alphabet $\mathcal{X}$ is *abstract* in general, it is understood that $g_n(\mathbf{x}) \equiv \frac{P_{X^n}(\mathbf{x})}{P_{\overline{X}^n}(\mathbf{x})}$ ($\mathbf{x} \in \mathcal{X}^n$) denotes the Radon-Nikodym derivative between two probability measures on $\mathcal{X}^n$ with values on a singular set assumed conventionally to be $+\infty$. Then, $\frac{P_{X^n}(X^n)}{P_{\overline{X}^n}(X^n)}$ is defined as $\frac{P_{X^n}(X^n)}{P_{\overline{X}^n}(X^n)} \equiv g_n(X^n)$, which is obviously a random variable. The probability distribution of the divergence-density rate is called the *divergence-spectrum* or more generally the *information-spectrum* (cf. Han and Verdú [1]).

[$\dagger$]For any sequence $\{Z_n\}_{n=1}^{\infty}$ of real-valued random variables, we define the *limit inferior in probability* (cf. Han and Verdú [1]) of $\{Z_n\}_{n=1}^{\infty}$ by $\text{p-}\liminf_{n \to \infty} Z_n = \sup\{\alpha | \lim_{n \to \infty} \Pr\{Z_n < \alpha\} = 0\}$.

*Proof:* If $R > \underline{D}(\mathbf{X}||\overline{\mathbf{X}})$, then by the definition of $\underline{D}(\mathbf{X}||\overline{\mathbf{X}})$ there exists an $0 < \varepsilon_0 < 1$ such that
$$\Pr\left\{\frac{1}{n}\log\frac{P_{X^n}(X^n)}{P_{\overline{X}^n}(X^n)} \leq R\right\} > \varepsilon_0$$
holds for infinitely many $n$. Hence,
$$\eta(R) \leq \liminf_{n\to\infty} \frac{1}{n}\log\frac{1}{\varepsilon_0} = 0.$$

$\square$

We now have the following quite general formula:

**Theorem 2.1** For any $r \geq 0$,
$$B_e(r|\mathbf{X}||\overline{\mathbf{X}}) = \inf_R \{R + \eta(R) \mid \eta(R) < r\}, \tag{2.2}$$
where $B_e(0|\mathbf{X}||\overline{\mathbf{X}}) = +\infty$ ($r = 0$).

**Remark 2.1** We notice here that $\eta(R) < r$ on the right-hand side of (2.2) is not $\eta(R) \leq r$. This is an essential difference, as will be seen in the proof below. Also, it is not difficult to check that $R + \eta(R) \geq 0$ for all $-\infty < R < +\infty$. $\square$

**Remark 2.2** Since it follows from Lemma 2.1 that
$$\inf_{R > \underline{D}(\mathbf{X}||\overline{\mathbf{X}})} \{R + \eta(R) \mid \eta(R) < r\} = \inf_{R > \underline{D}(\mathbf{X}||\overline{\mathbf{X}})} R$$
and inf on the right-hand side is attained by $R = \underline{D}(\mathbf{X}||\overline{\mathbf{X}})$, we may replace $\inf_R$ on the right-hand side of (2.2) by $\inf_{R \leq \underline{D}(\mathbf{X}||\overline{\mathbf{X}})}$ if $\eta(R)$ is continuous at $R = \underline{D}(\mathbf{X}||\overline{\mathbf{X}})$. $\square$

*Proof of Theorem 2.1*[‡]

) *Direct part:*

We use the notation that
$$S_n(a) = \left\{\mathbf{x} \in \mathcal{X}^n \,\bigg|\, \frac{1}{n}\log\frac{P_{X^n}(\mathbf{x})}{P_{\overline{X}^n}(\mathbf{x})} > a\right\}. \tag{2.3}$$
Let
$$\underline{R} = \inf\{R \mid \eta(R) < r\} \tag{2.4}$$

---
[‡]One of the referees suggested that the proof below based on the *information-spectrum slicing* is substantially similar to that of Varadhan's integral lemma (cf. Dembo and Zeitouni [4]), but this fact does *never* mean that Theorem 2.1 is a consequence of Varadhan's integral lemma, because the latter assumes the existence of a good rate function.

and consider the hypothesis testing with the acceptance region

$$\mathcal{A}_n = S_n(\underline{R} - \gamma)$$

with an arbitrarily small $\gamma > 0$. Then, the first kind of error probabbility is given by

$$\begin{aligned}\mu_n &= \Pr\{X^n \notin \mathcal{A}_n\} \\ &= \Pr\left\{\frac{1}{n}\log\frac{P_{X^n}(X^n)}{P_{\overline{X}^n}(X^n)} \leq \underline{R} - \gamma\right\}.\end{aligned}$$

Hence,

$$\liminf_{n\to\infty} \frac{1}{n}\log\frac{1}{\mu_n} = \eta(\underline{R} - \gamma).$$

On the other hand, (2.4) implies $\eta(\underline{R} - \gamma) \geq r$. Therefore,

$$\liminf_{n\to\infty} \frac{1}{n}\log\frac{1}{\mu_n} \geq r. \qquad (2.5)$$

Next, let us evaluate the second kind of error probability. First, put

$$\rho_0 = \inf_R \{R + \eta(R) \mid \eta(R) < r\}. \qquad (2.6)$$

We take $K$ large enough so as to satisfy $K > \rho_0$ and put $L = (K - \underline{R} + \gamma)/(2\gamma)$. Divide the interval $(\underline{R} - \gamma, K]$ into $L$ subintervals with equal width $2\gamma$ to define

$$I_i = (b_i - 2\gamma, b_i] \quad (i = 1, 2, \cdots, L), \qquad (2.7)$$

where $b_i \equiv \underline{R} - \gamma + 2i\gamma$. According to this interval partition, divide the set

$$T_0 = \left\{\mathbf{x} \in \mathcal{X}^n \,\middle|\, \underline{R} - \gamma < \frac{1}{n}\log\frac{P_{X^n}(\mathbf{x})}{P_{\overline{X}^n}(\mathbf{x})} \leq K\right\}$$

into the following $L$ subsets (*Information-spectrum slicing*):

$$S_n^{(i)} = \left\{\mathbf{x} \in \mathcal{X}^n \,\middle|\, \frac{1}{n}\log\frac{P_{X^n}(\mathbf{x})}{P_{\overline{X}^n}(\mathbf{x})} \in I_i\right\} \quad (i = 1, 2, \cdots, L).$$

Moreover, we define

$$S_n^{(0)} = \left\{\mathbf{x} \in \mathcal{X}^n \,\middle|\, \frac{1}{n}\log\frac{P_{X^n}(\mathbf{x})}{P_{\overline{X}^n}(\mathbf{x})} > K\right\}$$

to have

$$S_n(\underline{R} - \gamma) = \bigcup_{i=0}^{L} S_n^{(i)}. \qquad (2.8)$$

Since for $i = 1, 2, \cdots, L$ it holds that

$$\Pr\left\{X^n \in S_n^{(i)}\right\} \leq \Pr\left\{\frac{1}{n}\log\frac{P_{X^n}(X^n)}{P_{\overline{X}^n}(X^n)} \leq b_i\right\},$$

we have
$$\liminf_{n\to\infty} \frac{1}{n} \log \frac{1}{\Pr\left\{X^n \in S_n^{(i)}\right\}} \geq \eta(b_i).$$

Hence,
$$\Pr\left\{X^n \in S_n^{(i)}\right\} \leq e^{-n(\eta(b_i)-\gamma)} \quad (\forall n \geq n_0). \tag{2.9}$$

Moreover, if $\mathbf{x} \in S_n^{(i)}$ then
$$\frac{1}{n} \log \frac{P_{X^n}(\mathbf{x})}{P_{\overline{X}^n}(\mathbf{x})} > b_i - 2\gamma,$$

and so
$$P_{\overline{X}^n}(\mathbf{x}) \leq P_{X^n}(\mathbf{x}) e^{-n(b_i - 2\gamma)}.$$

As a result, by means of (2.9) we have §
$$\begin{aligned}\Pr\left\{\overline{X}^n \in S_n^{(i)}\right\} &\leq \sum_{\mathbf{x} \in S_n^{(i)}} P_{X^n}(\mathbf{x}) e^{-n(b_i - 2\gamma)} \\ &\leq e^{-n(b_i + \eta(b_i) - 3\gamma)}.\end{aligned} \tag{2.10}$$

Since $b_i \geq \underline{R} + \gamma$ for all $i = 1, 2, \cdots, L$,
$$b_i + \eta(b_i) \geq \rho_0 \quad (i = 1, 2, \cdots, L).$$

Substitution of this into (2.10) yields
$$\Pr\left\{\overline{X}^n \in S_n^{(i)}\right\} \leq e^{-n(\rho_0 - 3\gamma)} \quad (i = 1, 2, \cdots, L). \tag{2.11}$$

On the other hand, taking account that $\mathbf{x} \in S_n^{(0)}$ implies $P_{\overline{X}^n}(\mathbf{x}) \leq P_{X^n}(\mathbf{x}) e^{-nK}$, we have
$$\begin{aligned}\Pr\left\{\overline{X}^n \in S_n^{(0)}\right\} &= \sum_{\mathbf{x} \in S_n^{(0)}} P_{\overline{X}^n}(\mathbf{x}) \\ &\leq e^{-nK} \sum_{\mathbf{x} \in S_n^{(0)}} P_{X^n}(\mathbf{x}) \\ &\leq e^{-nK}.\end{aligned} \tag{2.12}$$

Consequently, from (2.8), (2.11), (2.12),
$$\lambda_n = \Pr\left\{\overline{X}^n \in S_n(\underline{R} - \gamma)\right\} \leq L e^{-n(\rho_0 - 3\gamma)} + e^{-nK}.$$

---

§In the case where the source alphabet $\mathcal{X}$ is *abstract* in general, the summation $\sum$ is understood to denote the integral $\int$.

We notice here that $K > \rho_0 - 3\gamma$ ($\gamma > 0$) because $K > \rho_0$. Thus,

$$\liminf_{n \to \infty} \frac{1}{n} \log \frac{1}{\lambda_n} \geq \rho_0 - 3\gamma,$$

which together with (2.5) concludes that $\rho_0 - 3\gamma$ is $r$-achievable (Notice here that $\gamma > 0$ is arbitrarily small).

) *Converse part:*¶

Let $\underline{R}$ and $\rho_0$ be defined as in (2.4), (2.6), respectively. Then, since $\eta(R)$ is monotone decreasing in $R$, there exists an $R_0$ such that $R_0 \geq \underline{R}$ and

$$\lim_{\varepsilon \downarrow 0}(R_0 + \varepsilon + \eta(R_0 + \varepsilon)) = \rho_0. \tag{2.13}$$

Let us consider the set

$$S_0 = \left\{ \mathbf{x} \in \mathcal{X}^n \,\bigg|\, \frac{1}{n} \log \frac{P_{X^n}(\mathbf{x})}{P_{\overline{X}^n}(\mathbf{x})} \leq R_0 + \gamma \right\},$$

where $\gamma > 0$ is an arbitrarily small constant. Then, by the definition of $\eta(R)$, there exists some divergent sequence $n_1 < n_2 < \cdots \to \infty$ of integers such that

$$\Pr\{X^{n_j} \in S_0\} \geq e^{-n_j(\eta(R_0+\gamma)+\tau)} \quad (\forall j \geq j_0), \tag{2.14}$$

where $\tau > 0$ is an arbitrarily small constant. Now let us use the contradiction argument. To do so, assume that $E = \rho_0 + 2\delta$ ($\delta > 0$ is a fixed constant) is $r$-achievable, i.e., assume that there exists an acceptance region $\mathcal{A}_n$ such that

$$\liminf_{n \to \infty} \frac{1}{n} \log \frac{1}{\mu_n} \geq r \tag{2.15}$$

and

$$\liminf_{n \to \infty} \frac{1}{n} \log \frac{1}{\lambda_n} \geq E \equiv \rho_0 + 2\delta. \tag{2.16}$$

Since $\mathbf{x} \in S_0$ implies

$$P_{X^n}(\mathbf{x}) \leq P_{\overline{X}^n}(\mathbf{x}) e^{n(R_0+\gamma)},$$

we have

$$\begin{aligned}
\Pr\{X^n \in S_0 \cap \mathcal{A}_n\} &= \sum_{\mathbf{x} \in S_0 \cap \mathcal{A}_n} P_{X^n}(\mathbf{x}) \\
&\leq \sum_{\mathbf{x} \in S_0 \cap \mathcal{A}_n} P_{\overline{X}^n}(\mathbf{x}) e^{n(R_0+\gamma)} \\
&\leq e^{n(R_0+\gamma)} \sum_{\mathbf{x} \in \mathcal{A}_n} P_{\overline{X}^n}(\mathbf{x}) \\
&= \lambda_n e^{n(R_0+\gamma)}.
\end{aligned} \tag{2.17}$$

---

¶Although it is usual in hypothesis testing problems to invoke the Neyman-Pearson lemma in order to prove the converse part, here we will give another simple elementary proof without recourse to the Neyman-Pearson lemma. This is to show that several alternative proofs are possible.

Furthermore, it follows from (2.16) that

$$\lambda_n \leq e^{-n(E-\gamma)} \quad (\forall n \geq n_0).$$

Substitution of this into (2.17) yields

$$\begin{aligned}\Pr\{X^n \in S_0 \cap \mathcal{A}_n\} &\leq e^{-n(E-R_0-2\gamma)} \\ &= e^{-n(\rho_0-R_0+2\delta-2\gamma)}.\end{aligned} \quad (2.18)$$

By virtue of (2.13), for any $\gamma > 0$ small enough,

$$\rho_0 \geq R_0 + \gamma + \eta(R_0 + \gamma) - \delta.$$

Therefore, by (2.18) we have

$$\Pr\{X^n \in S_0 \cap \mathcal{A}_n\} \leq e^{-n(\eta(R_0+\gamma)+\delta-\gamma)}.$$

Next, let us take $\tau > 0$, $\gamma > 0$ so small as to satisfy $\delta > 2\tau + \gamma$, then

$$\Pr\{X^n \in S_0 \cap \mathcal{A}_n\} \leq e^{-n(\eta(R_0+\gamma)+2\tau)}, \quad (2.19)$$

where $\tau > 0$ is the same one as in (2.14). On the other hand, by using (2.15), we obtain

$$\begin{aligned}\Pr\{X^n \in S_0 \cap \mathcal{A}_n^c\} &\leq \Pr\{X^n \in \mathcal{A}_n^c\} \\ &= \mu_n \leq e^{-n(r-\tau)} \quad (\forall n \geq n_0).\end{aligned} \quad (2.20)$$

We observe here that $\eta(R_0 + \gamma) < r$ for all $\gamma > 0$, and hence, for any sufficiently small $\tau > 0$,

$$\eta(R_0 + \gamma) + 2\tau < r - \tau.$$

Then, it follows from (2.19), (2.20) that

$$\begin{aligned}\Pr\{X^n \in S_0\} &= \Pr\{X^n \in S_0 \cap \mathcal{A}_n\} + \Pr\{X^n \in S_0 \cap \mathcal{A}_n^c\} \\ &\leq e^{-n(\eta(R_0+\gamma)+2\tau)} + e^{-n(r-\tau)} \\ &\leq 2e^{-n(\eta(R_0+\gamma)+2\tau)}\end{aligned} \quad (2.21)$$

for all $n \geq n_0$. However, since $\tau > 0$, (2.21) contradicts (2.14). Thus, the rate $E = \rho_0 + 2\delta$ cannot be $r$-achievable. Since $\delta > 0$ is arbitrary, it is concluded that any $E$ such that $E > \rho_0$ cannot be $r$-achievable. □

## 3 Examples

In this section we demonstrate several typical applications of Theorem 2.1. This is to verify the potentialities of Theorem 2.1.

**Example 3.1** Let the source alphabet $\mathcal{X}$ be *finite*, and consider the hypothesis testing where the null hypothsis $\mathbf{X} = (X_1, X_2, \cdots)$ and the alternative hypothesis $\overline{\mathbf{X}} = (\overline{X}_1, \overline{X}_2, \cdots)$ are stationary irreducible Markov sources subject to transition probabilities $P(x_2|x_1) = \Pr\{X_2 = x_2|X_1 = x_1\}$, $\overline{P}(x_2|x_1) = \Pr\{\overline{X}_2 = x_2|\overline{X}_1 = x_1\}$ $(x_1, x_2 \in \mathcal{X})$, respectively. Let $\mathcal{P}(\mathcal{X} \times \mathcal{X})$ denote the set of all probability distributions on $\mathcal{X} \times \mathcal{X}$, and, for any $Q \in \mathcal{P}(\mathcal{X} \times \mathcal{X})$ define the *conditional* divergences as

$$D(Q||P|q) = \sum_{x_1 \in \mathcal{X}} q(x_1) D(Q(\cdot|x_1)||P(\cdot|x_1)),$$

$$D(Q||\overline{P}|q) = \sum_{x_1 \in \mathcal{X}} q(x_1) D(Q(\cdot|x_1)||\overline{P}(\cdot|x_1)),$$

where $D(\cdot||\cdot)$ is the divergence (cf. Csiszár and Körner [6]), and $q(\cdot)$ and $Q(\cdot|\cdot)$ denote the marginal distribution and the conditional distribution of $Q$, respectively, which are defined as

$$q(x_1) = \sum_{x_2 \in \mathcal{X}} Q(x_1, x_2),$$

$$Q(x_2|x_1) = \frac{Q(x_1, x_2)}{q(x_1)}.$$

Then, by using Sanov theorem on the stationary irreducible Markov source (cf. Dembo and Zeitouni [4]), we have $\eta(R) = 0$ for $R \geq D(P||\overline{P}|p)$ ($p$ is the stationary distribution for $P$) and, for $R \leq D(P||\overline{P}|p)$,

$$\eta(R) = D(P_R||P|p_R), \tag{3.1}$$

$$R + \eta(R) = D(P_R||\overline{P}|p_R), \tag{3.2}$$

where, letting $\mathcal{P}_0$ be the set all probability distributions $Q \in \mathcal{P}(\mathcal{X} \times \mathcal{X})$ satisfying the *stationarity*, i.e.,

$$\mathcal{P}_0 = \left\{ Q \in \mathcal{P}(\mathcal{X} \times \mathcal{X}) \middle| \sum_{x_1 \in \mathcal{X}} Q(x_1, x) = \sum_{x_2 \in \mathcal{X}} Q(x, x_2) \text{ for all } x \in \mathcal{X} \right\}, \tag{3.3}$$

$P_R \in \mathcal{P}_0$ denotes the projection of $P$ on the plane:

$$\overline{\kappa}_R = \left\{ Q \in \mathcal{P}_0 \middle| \sum_{x_1, x_2 \in \mathcal{X}} Q(x_1, x_2) \log \frac{P(x_2|x_1)}{\overline{P}(x_2|x_1)} = R \right\} \tag{3.4}$$

as specified by

$$\inf_{Q \in \overline{\kappa}_R} D(Q||P|q) = D(P_R||P|p_R) \tag{3.5}$$

with $q$ being the marginal distribution of $Q$, and $p_R$ is the marginal distribution of $P_R$. Notice here that, since $Q$ moves on $\overline{\kappa}_R$, (3.5) implies also that

$$\inf_{Q \in \overline{\kappa}_R} D(Q||\overline{P}|q) = D(P_R||\overline{P}|p_R). \tag{3.6}$$

It is easy to see that $\underline{D}(\mathbf{X}||\overline{\mathbf{X}}) = D(P||\overline{P}|p)$ (cf. Barron [7]) and the function $\eta(R)$ given by (3.1) is continuous at $R = D(P||\overline{P}|p)$. Therefore, in view of Remark 2.2, it suffices to consider only $R$'s such that $R \le D(P||\overline{P}|p)$ on the right-hand side of (2.2). (Such an observation applies also to all the subsequent examples except for Example 3.6.) Thus, Theorem 2.1 leads us to

$$\begin{aligned} B_e(r|\mathbf{X}||\overline{\mathbf{X}}) &= \inf_R \left\{ D(P_R||\overline{P}|p_R) \mid D(P_R||P|p_R) < r \right\} \\ &= \inf_{Q \in \mathcal{P}_0 : D(Q||P|q) < r} D(Q||\overline{P}|q) \quad (\forall r > 0). \end{aligned} \qquad (3.7)$$

This result has been obtained by Natarajan [14]. This formula tells also that $B_e(r|\mathbf{X}||\overline{\mathbf{X}}) = 0$ whenever $r \ge D(\overline{P}||P|\overline{p})$ ($\overline{p}$ is the stationary distribution corresponding to $\overline{P}$).

If we consider the special case where sources $\mathbf{X}, \overline{\mathbf{X}}$ are both stationary memoryless subject to distributions $P, \overline{P}$ on $\mathcal{X}$, respectively, then formula (3.7) reduces to

$$B_e(r|\mathbf{X}||\overline{\mathbf{X}}) = \inf_{Q : D(Q||P) < r} D(Q||\overline{P}). \qquad (3.8)$$

This is nothing but Hoeffding's theorem [13] as is well known in the field of statistics. This tells also that $B_e(r|\mathbf{X}||\overline{\mathbf{X}}) = 0$ whenever $r \ge D(\overline{P}||P)$. □

**Example 3.2** Let us generalize Example 3.1 to the case with *unifilar* finite-state sources instead of stationary irreducible Markov sources. With the source alphabet $\mathcal{X}$ (finite) and the state set $\mathcal{S}$ (finite), let the null hypothesis $\mathbf{X} = \{X^n = (X_1, \cdot, X_n)\}_{n=1}^{\infty}$ be the unifilar finite-state source specified by

$$\begin{aligned} P_{X^n}(\mathbf{x}) &= \prod_{i=1}^n P(x_i|s_i) \quad (\mathbf{x} = (x_1, x_2, \cdots, x_n) \in \mathcal{X}^n) \qquad (3.9) \\ s_{i+1} &= f(x_i, s_i) \quad (s_i \in \mathcal{S}; \ i = 1, 2, \cdots, n, n+1), \qquad (3.10) \end{aligned}$$

and the let alternative hypothesis $\overline{\mathbf{X}} = \{\overline{X}^n = (\overline{X}_1, \cdot, \overline{X}_n)\}_{n=1}^{\infty}$ be the unifilar finite-state source specified by

$$\begin{aligned} P_{\overline{X}^n}(\mathbf{x}) &= \prod_{i=1}^n \overline{P}(x_i|s_i) \quad (\mathbf{x} = (x_1, x_2, \cdots, x_n) \in \mathcal{X}^n) \qquad (3.11) \\ s_{i+1} &= f(x_i, s_i) \quad (s_i \in \mathcal{S}; \ i = 1, 2, \cdots, n, n+1). \qquad (3.12) \end{aligned}$$

Given any fixed initial state $s_1 \in \mathcal{S}$, let $\mathcal{S}_0$ denote the set of all states $s \in \mathcal{S}$ that can be reached from $s_1$ with positive probability with respect to $P_{X^n}$. Next, letting $XS \equiv (X, S)$ be any random variable taking values in $\mathcal{X} \times \mathcal{S}_0$, put

$$S' = f(X, S). \qquad (3.13)$$

Moreover, let $\mathcal{V}_0$ denote the set of all the joint probability distributions $P_{XS}$ of random variables $XS$ satisfying both of the stationarity condition

$$P_{S'}(\cdot) = P_S(\cdot)$$

and the condition that the transition probability matrix $P_{S'|S}(\cdot|\cdot)$ is irreducible. Let the *projection* $P_{X_R S_R} \in \mathcal{V}_0$ of $P(\cdot|\cdot)$ on the plane $\overline{\lambda}_R$ be defined by

$$\inf_{P_{XS} \in \overline{\lambda}_R} D(P_{XS}||P|P_S) = D(P_{X_R S_R}||P|P_{S_R}), \tag{3.14}$$

where

$$\overline{\lambda}_R = \left\{ P_{XS} \in \mathcal{V}_0 \,\middle|\, \sum_{x \in \mathcal{X}, s \in \mathcal{S}_0} P_{XS}(x,s) \log \frac{P(x|s)}{\overline{P}(x|s)} = R \right\}. \tag{3.15}$$

Then, Sanov theorem on the unifilar finite-state source (cf. Han [2]) yields

$$\eta(R) = D(P_{X_R S_R}||P|P_{S_R}) \tag{3.16}$$
$$R + \eta(R) = D(P_{X_R S_R}||\overline{P}|P_{S_R}). \tag{3.17}$$

Notice here that, since $P_{XS}$ moves on $\overline{\lambda}_R$, (3.14) implies also that

$$\inf_{P_{XS} \in \overline{\lambda}_R} D(P_{XS}||\overline{P}|P_S) = D(P_{X_R S_R}||\overline{P}|P_{S_R}). \tag{3.18}$$

Thus, by Theorem 2.1 we have the following formula for the hypothesis testing $\mathbf{X}$ against $\overline{\mathbf{X}}$ with unifilar finite-state sources:

$$\begin{aligned}
B_e(r|\mathbf{X}||\overline{\mathbf{X}}) \\
&= \inf_R \left\{ D(P_{X_R S_R}||\overline{P}|P_{S_R}) \,|\, D(P_{X_R S_R}||P|P_{S_R}) < r \right\} \\
&= \inf_{P_{XS} \in \mathcal{V}_0 : D(P_{XS}||P|P_S) < r} D(P_{XS}||\overline{P}|P_S) \quad (\forall r > 0).
\end{aligned} \tag{3.19}$$

In the above argument we have taken account that in general the unifilar finite-state source is asymptotically a *mixture* of stationary or periodic *irreducible* sources. □

**Example 3.3** Let us consider the hypothesis testing with a *mixed* source as the null hypothesis, when the source alphabet $\mathcal{X}$ is *finite*. Let the alternative hypothesis $\overline{\mathbf{X}} = \{\overline{X}^n\}_{n=1}^\infty$ be a stationary memoryless source subject to probability distribution $\overline{P}$. Moreover, with any stationary memoryless sources $\mathbf{X}_1 = \{X_1^n\}_{n=1}^\infty, \mathbf{X}_2 = \{X_2^n\}_{n=1}^\infty$ subject to probability distributions $P_1, P_2$, respectively, let the null hypothesis $\mathbf{X} = \{X^n\}_{n=1}^\infty$ (called the *mixed source* of $\mathbf{X}_1$ and $\mathbf{X}_2$) be defined by

$$P_{X^n}(\mathbf{x}) = \alpha_1 P_{X_1^n}(\mathbf{x}) + \alpha_2 P_{X_2^n}(\mathbf{x}) \quad (\forall \mathbf{x} \in \mathcal{X}^n), \tag{3.20}$$

where $\alpha_1 > 0, \alpha_2 > 0$ are constants such that $\alpha_1 + \alpha_2 = 1$. In order to drive the required formula for this case, let the half-spaces $\nu_1, \nu_2$ be defined by

$$\nu_1 = \left\{ Q \in \mathcal{P}(\mathcal{X}) \,\middle|\, \sum_{x \in \mathcal{X}} Q(x) \log \frac{P_1(x)}{P_2(x)} \geq 0 \right\}, \tag{3.21}$$

$$\nu_2 = \left\{ Q \in \mathcal{P}(\mathcal{X}) \,\bigg|\, \sum_{x \in \mathcal{X}} Q(x) \log \frac{P_1(x)}{P_2(x)} \leq 0 \right\} \tag{3.22}$$

where $\mathcal{P}(\mathcal{X})$ is the set of all probability distributions on $\mathcal{X}$. Moreover, define other half-spaces in $\mathcal{P}(\mathcal{X})$ as

$$\kappa_R^{(1)} = \left\{ Q \in \mathcal{P}(\mathcal{X}) \,\bigg|\, \sum_{x \in \mathcal{X}} Q(x) \log \frac{P_1(x)}{\overline{P}(x)} \leq R \right\}, \tag{3.23}$$

$$\kappa_R^{(2)} = \left\{ Q \in \mathcal{P}(\mathcal{X}) \,\bigg|\, \sum_{x \in \mathcal{X}} Q(x) \log \frac{P_2(x)}{\overline{P}(x)} \leq R \right\}. \tag{3.24}$$

Then, letting the projections of $P_1, P_2$ on $\nu_1 \cap \kappa_R^{(1)}, \nu_2 \cap \kappa_R^{(2)}$ be denoted by $P_R^{(1)}, P_R^{(2)}$, respectively, Sanov theorem combined with the argument of *types* (cf. Han [20]) gives

$$\eta(R) = \min(D(P_R^{(1)} || P_1), D(P_R^{(2)} || P_2)). \tag{3.25}$$

Substituting this $\eta(R)$ into the right-hand side of (2.2) in Theorem 2.1, we can compute the value of $B_e(r|\mathbf{X}||\overline{\mathbf{X}})$ as a function of $r$ for the hypothesis testing with mixed sources.

Here, it easily follows from (3.25) that if $R \geq \min(D(P_1||\overline{P}), D(P_2||\overline{P}))$ then $\eta(R) = 0$, and that $\eta(R)$ is a monotone decreasing continuous function of $R$. Hence,

$$B_e(r|\mathbf{X}||\overline{\mathbf{X}}) \leq \min(D(P_1||\overline{P}), D(P_2||\overline{P})) \quad (\forall r > 0). \tag{3.26}$$

On the other hand, it follows again from (3.25) also that $\eta(h) > 0$ for any $h$ such that $h < \min(D(P_1||\overline{P}), D(P_2||\overline{P}))$, and so

$$\inf_R \{R + \eta(R) | \eta(R) < \eta(h)\} \geq h,$$

which implies that $h$ is $\eta(h)$-achievable. Hence, it holds that

$$\lim_{r \downarrow 0} B_e(r|\mathbf{X}||\overline{\mathbf{X}}) = \min(D(P_1||\overline{P}), D(P_2||\overline{P})). \tag{3.27}$$

□

**Remark 3.1** In fact, however, it is possible to drive a more general and much simpler formula for $B_e(r|\mathbf{X}||\overline{\mathbf{X}})$ with mixed sources, without any calculation of information spectra.. With *abstract* source alphabet $\mathcal{X}$ in general, let $\mathbf{X}_1 = \{X_1^n\}_{n=1}^\infty$, $\mathbf{X}_2 = \{X_2^n\}_{n=1}^\infty$, $\overline{\mathbf{X}}_1 = \{\overline{X}_1^n\}_{n=1}^\infty$, $\overline{\mathbf{X}}_2 = \{\overline{X}_2^n\}_{n=1}^\infty$ be any general sources. Consider the *mixed* source $\mathbf{X} = \{X^n\}_{n=1}^\infty$ of $\mathbf{X}_1$ and $\mathbf{X}_2$ and the *mixed* source $\overline{\mathbf{X}} = \{\overline{X}^n\}_{n=1}^\infty$ of $\overline{\mathbf{X}}_1$ and $\overline{\mathbf{X}}_2$, in the sense of Example 3.3, respectively. Then, for the hypothesis testing $\mathbf{X}$ against $\overline{\mathbf{X}}$, we have the general formula:

$$B_e(r|\mathbf{X}||\overline{\mathbf{X}}) = \min_{1 \leq i,j \leq 2} B_e(r|\mathbf{X}_i||\overline{\mathbf{X}}_j) \quad (\forall r > 0). \tag{3.28}$$

As for the detailed proof of (3.28), see Han [18]. □

**Example 3.4** Let us here consider the case with *countably infinite* source alphabet $\mathcal{X}$, say, $\mathcal{X} = \{1, 2, \cdots\}$. In this case, Sanov theorem as in Examples 3.1 ∼ 3.3 does not necessarily hold, while, since Cramér theorem (cf. Dembo and Zeitouni [4]) always holds, we can invoke here Cramér theorem instead of Sanov theorem. First, let $P = (p_1, p_2, \cdots)$, $\overline{P} = (\overline{p}_1, \overline{p}_2, \cdots)$ be any probablity distributions on $\mathcal{X}$, and let $X, \overline{X}$ denote the random variables such that $\Pr\{X = k\} = p_k$, $\Pr\{\overline{X} = k\} = \overline{p}_k$. Let $\mathbf{X} = \{X^n = (X_1, X_2, \cdots, X_n)\}_{n=1}^{\infty}$, $\overline{\mathbf{X}} = \{\overline{X}^n = (\overline{X}_1, \overline{X}_2, \cdots, \overline{X}_n)\}_{n=1}^{\infty}$ be the stationary memoryless sources specified by $X, \overline{X}$, respectively. Then, since the divergence-density rate is decomposed as

$$\frac{1}{n} \log \frac{P_{X^n}(X^n)}{P_{\overline{X}^n}(X^n)} = \frac{1}{n} \sum_{i=1}^{n} \log \frac{P_{X_i}(X_i)}{P_{\overline{X}_i}(X_i)}, \tag{3.29}$$

$\eta(R)$ in (2.1) can be expressed as

$$\eta(R) = \inf_{x \leq R} I(x), \tag{3.30}$$

where $I(x)$ is the large deviation rate function for (3.29). As usual, the moment generating function $M(\theta)$ of $\log \frac{P_X(X)}{P_{\overline{X}}(X)}$ is defined by

$$\begin{aligned} M(\theta) &= \mathrm{E} e^{\theta \log \frac{P_X(X)}{P_{\overline{X}}(X)}} = \sum_{i=1}^{\infty} p_i e^{\theta \log \frac{p_i}{\overline{p}_i}} \\ &= \sum_{i=1}^{\infty} p_i^{1+\theta} \overline{p}_i^{-\theta}. \end{aligned} \tag{3.31}$$

If we set $\Lambda(\theta) = \log M(\theta)$, Cramér theorem tells us that the rate function $I(x)$ is given by

$$I(x) = \sup_{\theta}(\theta x - \Lambda(\theta)), \tag{3.32}$$

where $-\log M(\theta)$ is called the Chernoff's $\theta$-distance (cf. Blahut [8], Cover and Thomas [9]). The expectation of $\log \frac{P_X(X)}{P_{\overline{X}}(X)}$ is computed as

$$\mathrm{E}\left[\log \frac{P_X(X)}{P_{\overline{X}}(X)}\right] = \sum_{i=1}^{\infty} p_i \log \frac{p_i}{\overline{p}_i} \equiv D(P||\overline{P}) \quad \text{(the divergence)}.$$

Thereofre, from (3.30) we see that if $R \geq D(P||\overline{P})$ then $\eta(R) = 0$, and if $R \leq D(P||\overline{P})$ then $\eta(R) = I(R)$. (It should be noted that $I(x)$ is monotone increasing in the range of $x \geq D(P||\overline{P})$, and monotone decreasing in the range of $x \leq D(P||\overline{P})$; and $I(x) = 0$ for $x = D(P||\overline{P})$.) Then, substituting (3.30) into (2.2) in Theorem 2.1, we can obtain the formula for computing the values of $B_e(r|\mathbf{X}||\overline{\mathbf{X}})$.

Substitution of (3.31) into (3.32) with $x = R$ yields

$$I(R) = \sup_{\theta}(\theta R - \log \sum_{i=1}^{\infty} p_i^{1+\theta} \overline{p}_i^{-\theta}), \tag{3.33}$$

which enables us to compute the values of $I(R)$. To compute this, differentiate the term in the bracket on the right-hand side of (3.33) with respect to $\theta$ and put it to zero to have the equation with respect to $\theta$:

$$R = \frac{\sum_{i=1}^{\infty} p_i^{1+\theta} \overline{p}_i^{-\theta} \log \frac{p_i}{\overline{p}_i}}{\sum_{i=1}^{\infty} p_i^{1+\theta} \overline{p}_i^{-\theta}} \equiv \varphi(\theta). \tag{3.34}$$

As far as $P \neq \overline{P}$, it is easy to check by using Schwarz inequality (cf. Gallager [10]) that $\varphi(\theta)$ on the right-hand side is continuous and strictly monotone increasing in $\theta$, because $M(\theta)$ is term-by-term continuously diffentiable (cf. Dembo and Zeitouni [4]). As a result, $\mathcal{D} \equiv \{-\infty < \varphi(\theta) < +\infty \mid \theta\}$ forms an interval on the real line. Therefore, in the case with $R \in \mathcal{D}$, $I(R)$ can be computed as

$$I(R) = \theta R - \log \sum_{i=1}^{\infty} p_i^{1+\theta} \overline{p}_i^{-\theta}, \tag{3.35}$$

where $\theta$ is the one as specified by (3.34). In this case, letting $\mathcal{P}(\mathcal{X})$ denote the set of all probability distributions on $\mathcal{X}$ and $Q_R$ denote the projection of the distribution $P$ on the plane $\kappa_R$ in $\mathcal{P}(\mathcal{X})$:

$$\kappa_R = \left\{ Q \in \mathcal{P}(\mathcal{X}) \,\middle|\, \sum_{i=1}^{\infty} Q(i) \log \frac{p_i}{\overline{p}_i} = R \right\},$$

we can ascertain by a direct calculation that

$$I(R) = D(Q_R \| P) \tag{3.36}$$

and

$$Q_R(i) = \frac{p_i^{1+\theta} \overline{p}_i^{-\theta}}{\sum_{i=1}^{\infty} p_i^{1+\theta} \overline{p}_i^{-\theta}} \quad (i \in \mathcal{X}) \tag{3.37}$$

with $\theta$ specified by the equation (3.34). Consequently, in the *countably infinite* alphabet case with $R \in \mathcal{D}$, Cramér theorem equivalently reduces to Sanov theorem as in (3.8) of Example 3.1 with *finite* alphabet. On the other hand, however, in the case with $R \notin \mathcal{D}$, the relation such as (3.36) does not hold. It then matters what interval $\mathcal{D}$ forms in general. In particular, if

$$D(\overline{P} \| P) < +\infty, \quad D(P \| \overline{P}) < +\infty, \tag{3.38}$$

then

$$[-D(\overline{P} \| P), D(P \| \overline{P})] \subset \mathcal{D}.$$

In this case, therefore, by using Sanov theorem in the same manner as in (3.8) of Example 3.1, we have for $0 < r \leq D(\overline{P} \| P)$ the formula

$$B_e(r | \mathbf{X} \| \overline{\mathbf{X}}) = \inf_{Q : D(Q \| P) < r} D(Q \| \overline{P}), \tag{3.39}$$

where it is easy to check that (3.39) holds also for $r > D(\overline{P}\|P)$ with $B_e(r|\mathbf{X}\|\overline{\mathbf{X}}) = 0$. The formula (3.39) gives an extended version with *countably infinite* source alphabet $\mathcal{X}$ of Hoeffding's theorem with *finite* source alphabet $\mathcal{X}$. It should be emphasized here that the formula (3.39) actually holds even with any *abstract* source alphabet $\mathcal{X}$ under the modest condition (3.38). In fact, the whole argument developed above continues to be valid, if only we equivalently rewrite $p_i^{1+\theta}\overline{p}_i^{-\theta}$ in the form $p_i\left(\frac{p_i}{\overline{p}_i}\right)^\theta$ where both of $\frac{p_i}{\overline{p}_i}$ and $\frac{\overline{p}_i}{p_i}$ in the latter form are well-defined as the Radon-Nikodym derivatives (cf. Billingsley [12]) even with any *abstract* source alphabet $\mathcal{X}$, in that the condition (3.38) is equivalent to the property that the probability measure $\overline{P}$ is *absolutely continuous* with respect to the probability measure $P$ and conversely the probability measure $P$ is *absolutely continuous* with respect to the probability measure $\overline{P}$.

The Cramér type equivalent of the formula (3.39) under condition (3.38) is found in Dembo and Zeitouni [4] where the Neyman-Pearson lemma is directly invoked, while here Theorem 2.1 is invoked. □

**Example 3.5** Let us consider the hypothesis testing where the null hypothesis $\mathbf{X}$ and the alternative hypothesis $\overline{\mathbf{X}}$ are both stationary memoryless sources subject to Gaussian distributions $N(\kappa, \sigma^2)$, $N(\overline{\kappa}, \sigma^2)$, respectively. Let the probability densities of these Gaussian distributions be written as

$$P_\kappa(x) = \frac{1}{\sqrt{2\pi}\sigma} e^{-\frac{(x-\kappa)^2}{2\sigma^2}},$$
$$P_{\overline{\kappa}}(x) = \frac{1}{\sqrt{2\pi}\sigma} e^{-\frac{(x-\overline{\kappa})^2}{2\sigma^2}}.$$

Denoting by $X$ the random variable subject to the probability density $P_\kappa$, the moment generating function $M(\theta) = \mathrm{E}(e^{\theta Y})$ of

$$Y = \log \frac{P_\kappa(X)}{P_{\overline{\kappa}}(X)} \tag{3.40}$$

is computed as

$$M(\theta) = e^{\frac{(\kappa-\overline{\kappa})^2(\theta+\theta^2)}{2\sigma^2}},$$

so that

$$\theta x - \log M(\theta) = \theta x - \frac{(\kappa-\overline{\kappa})^2(\theta+\theta^2)}{2\sigma^2}.$$

Then, the large deviation rate function $I(x)$ of (3.40) is given by

$$I(x) = \sup_\theta(\theta x - \log M(\theta)) = \frac{\sigma^2(x-a)^2}{2(\kappa-\overline{\kappa})^2}, \tag{3.41}$$

where, for simplicity, we have put $a = \frac{(\kappa - \overline{\kappa})^2}{2\sigma^2}$. Incidentally, we observe that $D(P_\kappa || P_{\overline{\kappa}}) = a$. Then, by means of Cramér theorem, $\eta(R)$ in Theorem 2.1 can be computed as

$$\begin{aligned}\eta(R) &= \inf_{x \leq R} I(x) = \min_{x \leq R} \frac{\sigma^2(x-a)^2}{2(\kappa - \overline{\kappa})^2} \\ &= \min\left\{[a-R]^+, \frac{\sigma^2(R-a)^2}{2(\kappa - \overline{\kappa})^2}\right\},\end{aligned} \quad (3.42)$$

from which it follows that

$$\begin{aligned}R + \eta(R) &= \min\left\{R + [a-R]^+, R + \frac{\sigma^2(R-a)^2}{2(\kappa - \overline{\kappa})^2}\right\} \\ &= \min\left\{R + [a-R]^+, \frac{\sigma^2(R+a)^2}{2(\kappa - \overline{\kappa})^2}\right\}.\end{aligned} \quad (3.43)$$

Thus, substituting (3.42) and (3.43) into the right-hand side of (2.2) in Theorem 2.1, we have

$$\begin{aligned}B_e(r|\mathbf{X}||\overline{\mathbf{X}}) &= \min\left\{[a-r]^+, (\sqrt{r} - \sqrt{a})^2\right\} \\ &= (\sqrt{r} - \sqrt{a})^2 \mathbf{1}[r \leq a],\end{aligned}$$

where $\mathbf{1}[\,\cdot\,]$ stands for the characteristic function. This formula tells us that $B_e(r|\mathbf{X}||\overline{\mathbf{X}})$ is monotone decreasing in $r$ for $0 < r < a$, and also that $B_e(0|\mathbf{X}||\overline{\mathbf{X}}) = a = D(P_\kappa || P_{\overline{\kappa}})$ and $B_e(r|\mathbf{X}||\overline{\mathbf{X}}) = 0$ for $r \geq a$. □

**Example 3.6** In all the examples that we have shown so far, the functions $\eta(R)$ were continuous in $R$. Here, we show an example in which $\eta(R)$ is discontinuous in $R$, where Remark 2.2 does not work. Let the source alphabet be $\mathcal{X} = \{0, 1\}$, and $S_n$ be a subset of $\mathcal{X}^n$ with size $|S_n| = 2^{\alpha n}$, where $\alpha$ is a constant such that $0 < \alpha < 1$. Moreover, let two elements $\mathbf{x}_0, \mathbf{x}_1 \in \mathcal{X}^n - S_n$ be fixed so that $\mathbf{x}_0 \neq \mathbf{x}_1$. The null hypothesis $\mathbf{X} = \{X^n\}_{n=1}^\infty$ be defined by

$$P_{X^n}(\mathbf{x}) = \begin{cases} 2^{-2\alpha n} & \text{for} & \mathbf{x} \in S_n, \\ 2^{-3\alpha n} & \text{for} & \mathbf{x} = \mathbf{x}_1, \\ 1 - 2^{-\alpha n} - 2^{-3\alpha n} & \text{for} & \mathbf{x} = \mathbf{x}_0, \\ 0 & \text{for} & \mathbf{x} \notin S_n \cup \{\mathbf{x}_1, \mathbf{x}_0\}, \end{cases} \quad (3.44)$$

where it is obvious that $P_{X^n}(S_n) = 2^{-\alpha n}$. The alternative hypothesis $\overline{\mathbf{X}} = \{\overline{X}^n\}_{n=1}^\infty$ be defined by $P_{\overline{X}^n}(\mathbf{x}) = 2^{-n}$ ($\forall \mathbf{x} \in \mathcal{X}^n$). Then, by a simple calculation, we see that the divergence-spectrum of this hypothesis testing consists of three points located at $1 + \frac{1}{n}\log(1 - 2^{-\alpha n} - 2^{-3\alpha n})$, $1 - 2\alpha$, $1 - 3\alpha$ with probabilities $1 - 2^{-\alpha n} - 2^{-3\alpha n}$, $2^{-\alpha n}, 2^{-3\alpha n}$, respectively. Therefore, by definition, the function $\eta(R)$ is given by

$$\eta(R) = \begin{cases} +\infty & \text{for} & R < 1 - 3\alpha, \\ 3\alpha & \text{for} & 1 - 3\alpha \leq R < 1 - 2\alpha, \\ \alpha & \text{for} & 1 - 2\alpha \leq R < 1, \\ 0 & \text{for} & 1 \leq R. \end{cases} \quad (3.45)$$

Hence, $R + \eta(R)$ is given by

$$R + \eta(R) = \begin{cases} +\infty & \text{for} & R < 1 - 3\alpha, \\ R + 3\alpha & \text{for} & 1 - 3\alpha \leq R < 1 - 2\alpha, \\ R + \alpha & \text{for} & 1 - 2\alpha \leq R < 1, \\ R & \text{for} & 1 \leq R. \end{cases} \quad (3.46)$$

Then, by Theorem 2.1, we have the formula

$$B_e(r|\mathbf{X}||\overline{\mathbf{X}}) = \begin{cases} 1 - \alpha & \text{for} & r > \alpha, \\ 1 & \text{for} & 0 < r \leq \alpha. \end{cases} \quad (3.47)$$

We observe here that, in the case of $r > \alpha$, $\inf_R$ on the right-hand side of (2.2) is attained by $R = R^\circ \equiv 1 - 2\alpha$, i.e.,

$$\inf_R \{R + \eta(R) \mid \eta(R) < r\} = R^\circ + \eta(R^\circ) \quad (R^\circ \equiv 1 - 2\alpha)$$
$$= 1 - \alpha.$$

In particular, we see that, if $r > 3\alpha$, $\inf_R$ is *not* attained by the boundary point $\underline{R} \equiv \inf\{R|\eta(R) < r\} = 1 - 3\alpha$ of $\{R|\eta(R) < r\}$, but by the internal point $R = R^\circ \equiv 1 - 2\alpha$. This kind of phenomenon has never taken place in the previous examples. Also, we should notice that formula (3.47) *cannot* be driven via the standard rate function method, differing from the previous examples, because in this case there does *not* exist any (lower semicontinuous) rate function. □

## 4 Hypothesis Testing and Large Deviation: Probability of Correct Testing

In this section we investigate the problem of determining the infimum $B_e^*(r|\mathbf{X}||\overline{\mathbf{X}})$ of achievable exponents for the second kind of *correct* probability $1 - \lambda_n$ under asymptotic constraints of the form $\mu_n \leq e^{-nr}$ on the first kind of error probability $\mu_n$ ($r > 0$ is a prescribed arbitrary constant), where $\lambda_n$ is the second kind of error probability. Let us first give the formal definititons, where $\mathbf{X} = \{X^n\}_{n=1}^\infty$, $\overline{\mathbf{X}} = \{\overline{X}^n\}_{n=1}^\infty$ indicate the null hypothesis and the alternative hypothesis, respectively.

**Definition 4.1** A rate $E$ is called $r$- *achievable* if there exists an acceptance region $\mathcal{A}_n$ such that

$$\liminf_{n \to \infty} \frac{1}{n} \log \frac{1}{\mu_n} \geq r \quad \text{and} \quad \limsup_{n \to \infty} \frac{1}{n} \log \frac{1}{1 - \lambda_n} \leq E.$$

**Definition 4.2 (The infimum of $r$-achievable correct exponents)**

$$B_e^*(r|\mathbf{X}||\overline{\mathbf{X}}) = \inf\{E \mid E \text{ is } r\text{-achievable}\}.$$

The purpose of this section is to determine $B_e^*(r|\mathbf{X}||\overline{\mathbf{X}})$ as a function of $r$. To this end, let us define the function $\eta(R)$ by

$$\eta(R) = \lim_{n \to \infty} \frac{1}{n} \log \frac{1}{\Pr\left\{\frac{1}{n} \log \frac{P_{X^n}(X^n)}{P_{\overline{X}^n}(X^n)} \leq R\right\}}. \tag{4.1}$$

This function is the same one as $\eta(R)$ defined by (2.1) in Section 2, but here we assume that the right-hand side of (4.1) has the limit. We notice here that $\eta(R)$ is monotone decreasing in $R$, and if $R > \underline{D}(\mathbf{X}||\overline{\mathbf{X}})$ then $\eta(R) = 0$ (cf. Lemma 2.1). The reason why we assume the existence of the limit in (4.1), on the contrary to in Section 2, will be made apparent below from the proof of Theorem 4.1.

Furthermore, for some technical reason, we assume in the sequel the following property about the information-spectrum that for any constant $M > 0$ there exists some constant $K > 0$ such that

$$\liminf_{n \to \infty} \frac{1}{n} \log \frac{1}{\Pr\left\{\frac{1}{n} \log \frac{P_{\overline{X}^n}(\overline{X}^n)}{P_{X^n}(\overline{X}^n)} \geq K\right\}} \geq M. \tag{4.2}$$

**Remark 4.1** This assumption[‖] means that the information-spectrum of $\overline{\mathbf{X}}$ with respect to $\mathbf{X}$ does not shift to the right faster than with any specified exponential speed of decay, when $n$ tends to $+\infty$. For example, if $\mathbf{X}, \overline{\mathbf{X}}$ are stationary memoryless sources with finite source alphabet subject to probability distributions $P_X, P_{\overline{X}}$, respectively, and there does not exist an $x \in \mathcal{X}$ for which $P_X(x) = 0$ and $P_{\overline{X}}(x) > 0$, then it is evident that the condition (4.2) is satisfied. □

We now have the following quite general formula, which is a *dual* counterpart of Theorem 2.1:

**Theorem 4.1** Assume that the limit in (4.1) exists and the condition (4.2) is satisfied. Then, for any $r \geq 0$,

$$B_e^*(r|\mathbf{X}||\overline{\mathbf{X}}) = \inf_R \left\{R + \eta(R) + [r - \eta(R)]^+\right\}, \tag{4.3}$$

where $[x]^+ = \max(x, 0)$ and we have put $B_e^*(0|\mathbf{X}||\overline{\mathbf{X}}) = 0$ ($r = 0$).

**Remark 4.2** Since it follows from Lemma 2.1 that

$$\inf_{R > \underline{D}(\mathbf{X}||\overline{\mathbf{X}})} \left\{R + \eta(R) + [r - \eta(R)]^+\right\} = \inf_{R > \underline{D}(\mathbf{X}||\overline{\mathbf{X}})} (R + r),$$

the inf on the right-hand side is attained by $R = \underline{D}(\mathbf{X}||\overline{\mathbf{X}})$. Therefore, $\inf_R$ on the right-hand side of (4.3) may be replaced by $\inf_{R \leq \underline{D}(\mathbf{X}||\overline{\mathbf{X}})}$ if $\eta(R)$ is continuous at $R = \underline{D}(\mathbf{X}||\overline{\mathbf{X}})$. □

---

[‖]One of the referees suggests the striking similarity between the condition (4.2) and the standard concept of *exponential tightness* in large deviation theory (e.g., cf. Dembo and Zeitouni [4]).

*Proof of Theorem 4.1.*

) *Direct part:*

In the proof of the direct part we do not need the assumtion (4.2). First, keep in mind that $\eta(R)$ in
$$R + \eta(R) + [r - \eta(R)]^+$$
on the right-hand side of (4.3) is monotone decreasing, and set
$$\rho_0^* = \inf_R \left\{ R + \eta(R) + [r - \eta(R)]^+ \right\}. \tag{4.4}$$

Then, there exists an $R_0$ such that $\rho_0^*$ is expressed as
$$\rho_0^* = \lim_{\varepsilon \downarrow 0}(R_0 + \varepsilon + \eta(R_0 + \varepsilon) + [r - \eta(R_0 + \varepsilon)]^+), \tag{4.5}$$

which we rewrite as
$$\rho_0^* = R_0 + \gamma + \eta(R_0 + \gamma) + [r - \eta(R_0 + \gamma)]^+ - \nu(\gamma), \tag{4.6}$$

where $\gamma > 0$ is an arbitrarily small constant and $\nu(\gamma) \to 0$ as $\gamma \to 0$. We use here the notation that
$$S_n^*(a) = \left\{ \mathbf{x} \in \mathcal{X}^n \,\middle|\, \frac{1}{n} \log \frac{P_{X^n}(\mathbf{x})}{P_{\overline{X}^n}(\mathbf{x})} \leq a \right\}. \tag{4.7}$$

Then, since the existence of the limit in (4.1) was assumed, we have
$$e^{-n(\eta(R_0+\gamma)+\tau)} \leq \Pr\{X^n \in S_n^*(R_0 + \gamma)\} \leq e^{-n(\eta(R_0+\gamma)-\tau)} \quad (\forall n \geq n_0), \tag{4.8}$$

where $\tau > 0$ is an arbitrarily small constant. Next, define a subset $\mathcal{C}_n$ of $S_n^*(R_0 + \gamma)$ as follows; if $\eta(R_0 + \gamma) \geq r$ then set $\mathcal{C}_n = S_n^*(R_0 + \gamma)$, otherwise if $\eta(R_0 + \gamma) < r$ then set $\mathcal{C}_n = T_n$ where $T_n$ is any subset of $S_n^*(R_0 + \gamma)$ such that
$$\lim_{n \to \infty} \frac{1}{n} \log \frac{1}{\Pr\{X^n \in T_n\}} = r. \tag{4.9}$$

It should be noted here that it is always possible to choose such a subset $T_n$, because in the case with $\eta(R_0 + \gamma) < r$ we can make $\eta(R_0 + \gamma) + \tau < r$ hold with $\tau > 0$ small enough, where we may consider a *randomized* hypothesis testing if necessary. Now, consider the hypothesis testing with $\mathcal{C}_n$ as the critical region. First, we evaluate the value of the first kind of error probablity $\mu_n$. In the case with $\eta(R_0 + \gamma) \geq r$, since $\mathcal{C}_n = S_n^*(R_0 + \gamma)$, by means of (4.8) we have
$$\begin{aligned}\Pr\{X^n \in \mathcal{C}_n\} &\leq e^{-n(\eta(R_0+\gamma)-\tau)} \\ &\leq e^{-n(r-\tau)} \quad (\forall n \geq n_0),\end{aligned}$$

while in the case with $\eta(R_0 + \gamma) < r$, by means of (4.9) we have
$$\Pr\{X^n \in \mathcal{C}_n\} \leq e^{-n(r-\tau)} \quad (\forall n \geq n_0).$$

Then, in either case, it holds that

$$\Pr\{X^n \in \mathcal{C}_n\} \leq e^{-n(r-\tau)}. \tag{4.10}$$

Therefore, the first kind of error probablity $\mu_n$ is evaluated as

$$\mu_n \equiv \Pr\{X^n \in \mathcal{C}_n\} \leq e^{-n(r-\tau)}.$$

Hence,

$$\liminf_{n \to \infty} \frac{1}{n} \log \frac{1}{\mu_n} \geq r - \tau.$$

Since $\tau > 0$ is arbitrary, it is concluded that

$$\liminf_{n \to \infty} \frac{1}{n} \log \frac{1}{\mu_n} \geq r. \tag{4.11}$$

Next, we evaluate the value of the second kind of correct probability $1 - \lambda_n$, where $\lambda_n$ is the second kind of error probability. First, we observe that if $\mathbf{x} \in S_n^*(R_0 + \gamma)$ then

$$P_{\overline{X}^n}(\mathbf{x}) \geq P_{X^n}(\mathbf{x}) e^{-n(R_0+\gamma)} \tag{4.12}$$

holds. Then, in the case with $\eta(R_0 + \gamma) \geq r$, since $\mathcal{C}_n = S_n^*(R_0 + \gamma)$, it follows from (4.8) that

$$\begin{aligned}
\Pr\left\{\overline{X}^n \in \mathcal{C}_n\right\} &= \sum_{\mathbf{x} \in \mathcal{C}_n} P_{\overline{X}^n}(\mathbf{x}) \\
&\geq \sum_{\mathbf{x} \in \mathcal{C}_n} P_{X^n}(\mathbf{x}) e^{-n(R_0+\gamma)} \\
&= e^{-n(R_0+\gamma)} \Pr\{X^n \in S_n^*(R_0 + \gamma)\} \\
&\geq e^{-n(R_0+\gamma+\eta(R_0+\gamma)+\tau)} \quad (\forall n \geq n_0).
\end{aligned} \tag{4.13}$$

Similarly, in the case with $\eta(R_0 + \gamma) < r$, since $\mathcal{C}_n = T_n$, it follows from (4.9) that

$$\Pr\left\{\overline{X}^n \in \mathcal{C}_n\right\} \geq e^{-n(R_0+\gamma+r+\tau)} \quad (\forall n \geq n_0). \tag{4.14}$$

Summarizing (4.13) and (4.14), in either case we have

$$\Pr\left\{\overline{X}^n \in \mathcal{C}_n\right\} \geq e^{-n(R_0+\gamma+\eta(R_0+\gamma)+[r-\eta(R_0+\gamma)]^+ + \tau)}. \tag{4.15}$$

Substitution of (4.6) into (4.15) yields

$$\Pr\left\{\overline{X}^n \in \mathcal{C}_n\right\} \geq e^{-n(\rho_0^* + \tau + \nu(\gamma))}.$$

Hence,

$$\begin{aligned}
1 - \lambda_n &= \Pr\left\{\overline{X}^n \in \mathcal{C}_n\right\} \\
&\geq e^{-n(\rho_0^* + \tau + \nu(\gamma))},
\end{aligned}$$

from which it follows that
$$\limsup_{n\to\infty} \frac{1}{n}\log\frac{1}{1-\lambda_n} \le \rho_0^* + \tau + \nu(\gamma). \tag{4.16}$$
We notice here that we can make $\tau + \nu(\gamma) \to 0$, because $\tau > 0$ and $\gamma > 0$ are both made arbitrarily small. Thus, by virtue of (4.11) and (4.16) we conclude that any rate $E$ such that $E > \rho_0^*$ is $r$-ahievable.

) *Converse part:*

In the proof of the converse part we need the assumption (4.2). First, let $K > 0$ be a constant large enough (to be specified below) and $\gamma > 0$ be an arbitrarily small constant. Putting $L = \frac{2K}{\gamma}$, we divide the interval $(-K, K]$ into $L$ subintervals with equal width $\gamma$ to have
$$I_i = (c_i - \gamma, c_i] \quad (i = 1, 2, \cdots, L),$$
where $c_i \equiv K - (i-1)\gamma$. According to this interval partition, divide the set
$$T_n^* = \left\{ \mathbf{x} \in \mathcal{X}^n \left| -K < \frac{1}{n}\log\frac{P_{X^n}(\mathbf{x})}{P_{\overline{X}^n}(\mathbf{x})} \le K \right. \right\}$$
into the $L$ subsets
$$S_n^{(i)} = \left\{ \mathbf{x} \in \mathcal{X}^n \left| \frac{1}{n}\log\frac{P_{X^n}(\mathbf{x})}{P_{\overline{X}^n}(\mathbf{x})} \in I_i \right. \right\} \quad (i = 1, 2, \cdots, L).$$
This operation is called the *information-spectrum slicing*. Moreover, we define
$$S_n^{(0)} = \left\{ \mathbf{x} \in \mathcal{X}^n \left| \frac{1}{n}\log\frac{P_{X^n}(\mathbf{x})}{P_{\overline{X}^n}(\mathbf{x})} \le -K \right. \right\},$$
$$S_n^{(-1)} = \left\{ \mathbf{x} \in \mathcal{X}^n \left| \frac{1}{n}\log\frac{P_{X^n}(\mathbf{x})}{P_{\overline{X}^n}(\mathbf{x})} > K \right. \right\},$$
where it is obvious that $\mathcal{X}^n = \bigcup_{j=-1}^{L} S_n^{(j)}$. Suppose that $E$ is $r$-achievable, i.e., that there exists a critical region $\mathcal{C}_n$ such that
$$\liminf_{n\to\infty} \frac{1}{n}\log\frac{1}{\mu_n} \ge r, \tag{4.17}$$
$$\limsup_{n\to\infty} \frac{1}{n}\log\frac{1}{1-\lambda_n} \le E. \tag{4.18}$$
Then, from (4.17) we have
$$\mu_n \le e^{-n(r-\tau)} \quad (\forall n \ge n_0), \tag{4.19}$$
where $\tau > 0$ is an arbitrarily small constant. In order to evaluate the value of $\Pr\left\{\overline{X}^n \in \mathcal{C}_n\right\}$, let us first evaluate the value of
$$\Pr\left\{X^n \in \mathcal{C}_n^{(i)}\right\} \quad (i = 1, 2, \cdots, L),$$

where $\mathcal{C}_n^{(i)} \equiv S_n^{(i)} \cap \mathcal{C}_n$ ($i = -1, 0, 1, 2, \cdots, L$). We now evaluate the value of $\Pr\left\{X^n \in \mathcal{C}_n^{(i)}\right\}$ ($i = 1, 2, \cdots, L$) in two ways as follows. First, we observe that

$$\Pr\left\{X^n \in \mathcal{C}_n^{(i)}\right\} \leq \Pr\{X^n \in \mathcal{C}_n\} = \mu_n,$$

which together with (4.19) yields

$$\Pr\left\{X^n \in \mathcal{C}_n^{(i)}\right\} \leq e^{-n(r-\tau)}. \tag{4.20}$$

Next, by the definitions of $\eta(c_i)$ and $S_n^{(i)}$, we see that

$$\begin{aligned}
\Pr\{X^n \in S_n^{(i)}\} &\leq \Pr\left\{\frac{1}{n}\log\frac{P_{X^n}(X^n)}{P_{\overline{X}^n}(X^n)} \leq c_i\right\} \\
&\leq e^{-n(\eta(c_i)-\tau)} \quad (\forall n \geq n_0).
\end{aligned}$$

Hence,

$$\begin{aligned}
\Pr\left\{X^n \in \mathcal{C}_n^{(i)}\right\} &\leq \Pr\left\{X^n \in S_n^{(i)}\right\} \\
&\leq e^{-n(\eta(c_i)-\tau)}. 
\end{aligned} \tag{4.21}$$

A consequence of (4.20) and (4.21) is

$$\Pr\left\{X^n \in \mathcal{C}_n^{(i)}\right\} \leq e^{-n(\eta(c_i)+[r-\eta(c_i)]^+ -\tau)} \quad (i = 1, 2, \cdots, L). \tag{4.22}$$

We can now evaluate the value of $\Pr\left\{\overline{X}^n \in \mathcal{C}_n^{(i)}\right\}$ as follows. Since $\mathbf{x} \in \mathcal{C}_n^{(i)}$ implies $\mathbf{x} \in S_n^{(i)}$ ($i = 1, 2, \cdots, L$) and hence also $P_{\overline{X}^n}(\mathbf{x}) \leq P_{X^n}(\mathbf{x})e^{-n(c_i-\gamma)}$, we have

$$\begin{aligned}
\Pr\left\{\overline{X}^n \in \mathcal{C}_n^{(i)}\right\} &= \sum_{\mathbf{x} \in \mathcal{C}_n^{(i)}} P_{\overline{X}^n}(\mathbf{x}) \\
&\leq \sum_{\mathbf{x} \in \mathcal{C}_n^{(i)}} P_{X^n}(\mathbf{x})e^{-n(c_i-\gamma)} \\
&= e^{-n(c_i-\gamma)}\Pr\left\{X^n \in \mathcal{C}_n^{(i)}\right\} \\
&\leq e^{-n(c_i+\eta(c_i)+[r-\eta(c_i)]^+ -\gamma-\tau)} 
\end{aligned} \tag{4.23}$$

for $i = 1, 2, \cdots, L$, where we have used (4.22) in the last inequality. Furthermore, let us evaluate the values of $\Pr\left\{\overline{X}^n \in S_n^{(-1)}\right\}$ and $\Pr\left\{\overline{X}^n \in S_n^{(0)}\right\}$. Since $\mathbf{x} \in S_n^{(-1)}$ implies $P_{\overline{X}^n}(\mathbf{x}) \leq P_{X^n}(\mathbf{x})e^{-nK}$, we obtain

$$\begin{aligned}
\Pr\left\{\overline{X}^n \in S_n^{(-1)}\right\} &= \sum_{\mathbf{x} \in S_n^{(-1)}} P_{\overline{X}^n}(\mathbf{x}) \\
&\leq \sum_{\mathbf{x} \in S_n^{(-1)}} P_{X^n}(\mathbf{x})e^{-nK} \\
&\leq e^{-nK}. 
\end{aligned} \tag{4.24}$$

Recalling here that

$$\Pr\left\{\overline{X}^n \in S_n^{(0)}\right\} = \Pr\left\{\frac{1}{n}\log\frac{P_{X^n}(\overline{X}^n)}{P_{\overline{X}^n}(\overline{X}^n)} \leq -K\right\}$$

$$= \Pr\left\{\frac{1}{n}\log\frac{P_{\overline{X}^n}(\overline{X}^n)}{P_{X^n}(\overline{X}^n)} \geq K\right\}.$$

and noting the assumption (4.2), we see that for any $M > 0$ there exists a $K > 0$ large enough such that

$$\Pr\left\{\overline{X}^n \in S_n^{(0)}\right\} \leq e^{-n(M-\tau)} \quad (\forall n \geq n_0). \tag{4.25}$$

Summarizing up (4.23)~(4.25), we have

$$1 - \lambda_n$$
$$= \Pr\left\{\overline{X}^n \in \mathcal{C}_n\right\} = \sum_{i=-1}^{L} \Pr\left\{\overline{X}^n \in \mathcal{C}_n^{(i)}\right\}$$
$$\leq \sum_{i=1}^{L} e^{-n(c_i + \eta(c_i) + [r - \eta(c_i)]^+ - \gamma - \tau)} + e^{-nK} + e^{-n(M-\tau)}. \tag{4.26}$$

On the other hand, since, by the definition (4.4) of $\rho_0^*$,

$$c_i + \eta(c_i) + [r - \eta(c_i)]^+ \geq \rho_0^* \quad (i = 1, 2, \cdots, L),$$

it follows from (4.26) that

$$1 - \lambda_n \leq L e^{-n(\rho_0^* - \gamma - \tau)} + e^{-nK} + e^{-n(M-\tau)}.$$

Thus, if we take $M > 0$ and $K > 0$ large enough, then

$$\limsup_{n \to \infty} \frac{1}{n}\log\frac{1}{1 - \lambda_n} \geq \rho_0^* - \gamma - \tau. \tag{4.27}$$

Therefore, $E \geq \rho_0^* - \gamma - \tau$ holds, owing to (4.18), (4.27). Since both of $\gamma > 0$ and $\tau > 0$ are arbitrary, we can let $\gamma \to 0$, $\tau \to 0$ to get $E \geq \rho_0^*$. Thus, it is concluded that any $r$-achievable rate $E$ cannot be smaller than $\rho_0^*$. □

## 5 Examples

In this section we demonstrate several typical applications of Theorem 4.1. This is to verify the potentialities of Theorem 4.1.

**Example 5.1** As in Example 3.1, let us consider the hypothesis testing with stationary irreducible Markov sources $\mathbf{X}, \overline{\mathbf{X}}$ with *finite* source alphabet. With the same notation as in Example 3.1, it follows also here with Sanov theorem that (3.1) and (3.2) hold, i.e., $\eta(R) = 0$ for $R \geq D(P||\overline{P}|p)$ and, for $R \leq D(P||\overline{P}|p)$,

$$\eta(R) = D(P_R||P|p_R),$$
$$R + \eta(R) = D(P_R||\overline{P}|p_R),$$

and so, by Theorem 4.1 we have

$$B_e^*(r|\mathbf{X}||\overline{\mathbf{X}}) = \inf_R \left\{ D(P_R||\overline{P}|p_R) + [r - D(P_R||P|p_R)]^+ \right\}. \tag{5.1}$$

It is easy to check that, if $r < D(\overline{P}||P|\overline{p})$ ($\overline{p}$ is the stationary distriibution for $\overline{P}$) then $B_e^*(r|\mathbf{X}||\overline{\mathbf{X}}) = 0$, whereas if $r \geq D(\overline{P}||P|\overline{p})$ then $\inf_R$ on the right-hand side of (5.1) is attained by an $R$ such that

$$D(P_R||P|p_R) \leq r,$$

and hence in this latter case we have

$$B_e^*(r|\mathbf{X}||\overline{\mathbf{X}}) = \inf_{R:D(P_R||P|p_R)\leq r} \left\{ D(P_R||\overline{P}|p_R) + r - D(P_R||P|p_R) \right\}$$
$$= \inf_{Q\in\mathcal{P}_0:D(Q||P|q)\leq r} \left\{ D(Q||\overline{P}|q) + r - D(Q||P|q) \right\}. \tag{5.2}$$

This formula has been driven by Nakagawa and Kanaya [16]. Here, $B_e^*(r|\mathbf{X}||\overline{\mathbf{X}}) = 0$ whenever $r \leq D(\overline{P}||P|\overline{p})$.

Let us consider the special case where sources $\mathbf{X}, \overline{\mathbf{X}}$ are both stationary memoryless subject to probability distributions $P, \overline{P}$ on $\mathcal{X}$, respectively. Then, in the case of $r \geq D(\overline{P}||P)$, formula (5.2) reduces to

$$B_e^*(r|\mathbf{X}||\overline{\mathbf{X}}) = \inf_{Q:D(Q||P)\leq r} \left\{ D(Q||\overline{P}) + r - D(Q||P) \right\}. \tag{5.3}$$

This formula has first been established by Han and Kobayashi [15] based on the method of types. On the other hand, we have $B_e^*(r|\mathbf{X}||\overline{\mathbf{X}}) = 0$ whenever $r < D(\overline{P}||P)$.
□

**Example 5.2** In order to generalize Example 5.1, as in Example 3.2 of §2 let us consider the hypothesis testing with *unifilar* finite-state sources $\mathbf{X}, \overline{\mathbf{X}}$ as the null and alternative hypotheses, respectively. With the same notation as in Example 3.2, Theorem 4.1 together with (3.16) and (3.17) gives the formula for the hypothesis testing $\mathbf{X}$ against $\overline{\mathbf{X}}$:

$$B_e^*(r|\mathbf{X}||\overline{\mathbf{X}})$$
$$= \inf_R \left\{ D(P_{X_R S_R}||\overline{P}|P_{S_R}) + [r - D(P_{X_R S_R}||P|P_{S_R})]^+ \right\}$$
$$= \inf_{P_{XS}\in\mathcal{V}_0} \left\{ D(P_{XS}||\overline{P}|P_S) + [r - D(P_{XS}||P|P_S)]^+ \right\}. \tag{5.4}$$

**Example 5.3** Let the source alphabet $\mathcal{X}$ be *finite*, and, as in Example 3.3, let us consider the hypothesis testing with a *mixed* source $\mathbf{X} = \{X^n\}_{n=1}^{\infty}$ as the null hypothesis and a stationary memoryless source $\overline{\mathbf{X}} = \{\overline{X}^n\}_{n=1}^{\infty}$ subject to probability distribution $\overline{P}$ as the alternative hypothesis. In order to satisfy the assumption (4.2), let $P_1(x) > 0$, $P_2(x) > 0$ ($\forall x \in \mathcal{X}$). Here, recall that the mixed source $\mathbf{X} = \{X^n\}_{n=1}^{\infty}$ was defined as

$$P_{X^n}(\mathbf{x}) = \alpha_1 P_{X_1^n}(\mathbf{x}) + \alpha_2 P_{X_2^n}(\mathbf{x}) \quad (\forall \mathbf{x} \in \mathcal{X}^n), \tag{5.5}$$

where $\mathbf{X}_1 = \{X_1^n\}_{n=1}^{\infty}, \mathbf{X}_2 = \{X_2^n\}_{n=1}^{\infty}$ are stationary memoryless sources subject to probability distributions $P_1, P_2$, respectively. Define $\nu_1, \nu_2, \kappa_R^{(1)}, \kappa_R^{(2)}$ as in (3.21)$\sim$ (3.24) of Example 3.3, and similarly, let the projections of $P_1, P_2$ on $\nu_1 \cap \kappa_R^{(1)}, \nu_2 \cap \kappa_R^{(2)}$ be denoted by $P_R^{(1)}, P_R^{(2)}$, respectively. Then, application of Sanov theorem gives

$$\eta(R) = \min(D(P_R^{(1)}||P_1), D(P_R^{(2)}||P_2)), \tag{5.6}$$

from which we see that if $R \geq \min(D(P_1||\overline{P}), D(P_2||\overline{P}))$ then $\eta(R) = 0$.

Finally, by substituting $\eta(R)$ of (5.6) into the right-hand side of (4.3) in Theorem 4.1 we have the computable formula for $B_e^*(r|\mathbf{X}||\overline{\mathbf{X}})$ as a function of $r$.  □

**Remark 5.1** Unfortunately, for $B_e^*(r|\mathbf{X}||\overline{\mathbf{X}})$ that we are considering here, such a simple formula for mixed sources as (3.28) in Remark 3.1 does not hold.  □

**Remark 5.2** So far, we have considered only the case with *finite* source alphabet $\mathcal{X}$ where Sanov theorem played the key role. On the other hand, in the case of general stationary memoryless sources with *countably infinite* or *abstract* source alphabet $\mathcal{X}$, Sanov theorem does not necessarily hold. However, since Cramér theorem always works, we can invoke Cramér theorem, instead of Sanov theorem, in order to compute the value of $B_e^*(r|\mathbf{X}||\overline{\mathbf{X}})$, when $\mathbf{X}, \overline{\mathbf{X}}$ are both stationary memoryless sources. Then, it suffices to use the same rate function $I(x)$ as specified in (3.30) of Example 3.4, i.e.,

$$\eta(R) = \inf_{x \leq R} I(x). \tag{5.7}$$

With the same notation as in Example 3.4, we see that we can write the right-hand side of (5.7) in terms of divergences (with Sanov theorem) only when $R \in \mathcal{D}$.  □

**Example 5.4** Let us consider the hypothesis testing with stationary memoryless Gaussian sources $\mathbf{X} = \{P_\kappa\}, \overline{\mathbf{X}} = \{P_{\overline{\kappa}}\}$ as in Example 3.5. Since $\eta(R)$ and $R + \eta(R)$ are given by (3.42), (3.43), substitution of these (3.42), (3.43) into (4.3) in Theorem 4.1 and some simple calculation yield the formula

$$B_e^*(r|\mathbf{X}||\overline{\mathbf{X}}) = (\sqrt{r} - \sqrt{a})^2 \mathbf{1}[r \geq a], \tag{5.8}$$

where $a = D(P_\kappa || P_{\overline{\kappa}})$. We noitce here that the function (5.8) is symmetric to the function $B_e(r|\mathbf{X}||\overline{\mathbf{X}})$ in Example 3.5 with respect to the $y$-axis. The formula (5.8) tells us that $B_e^*(r|\mathbf{X}||\overline{\mathbf{X}})$ is a monotone increasing function of $r$, and that $B_e^*(r|\mathbf{X}||\overline{\mathbf{X}}) = 0$ whenever $r \leq a$. □

# 6 Generalized Hypothesis Testing

So far, we have studied the hypothesis testing problem with general sources $\mathbf{X} = \{X^n\}_{n=1}^\infty$, $\overline{\mathbf{X}} = \{\overline{X}^n\}_{n=1}^\infty$ as null and alternative hypotheses, respectively. However, it is easy to observe that Theorem 2.1 and Theorem 4.1 in the previous sections continue to be valid as they are, even if we replace the probability distribution $P_{\overline{X}^n}$ of the *alternative* hypothesis by any *nonnegative measure* $G_n$ with $G_n(\emptyset) = 0$ (*not necessarily* a probability measure), where the second kind of error probability $\lambda_n \equiv \Pr\{\overline{X}^n \in \mathcal{A}_n\}$ should be interpreted in turn as denoting the value of the nonnegative measure $\lambda_n \equiv G_n(\mathcal{A}_n)$. This is called the *generalized* hypothesis testing. Then, if we define

$$\kappa_n \equiv G_n(\mathcal{X}^n) \quad (n = 1, 2, \cdots),$$
$$\underline{\kappa} \equiv \limsup_{n \to \infty} \frac{1}{n} \log \frac{1}{G_n(\mathcal{X}^n)},$$

Theorem 4.1 is meaningful only when $\underline{\kappa} < +\infty$, where $B_e^*(0|\mathbf{X}||\overline{\mathbf{X}}) = 0$ in Theorem 4.1 needs to be replaced by $B_e^*(0|\mathbf{X}||\overline{\mathbf{X}}) = \underline{\kappa}$, and $1 - \lambda_n$ in Definition 4.1 needs to be replaced by $\kappa_n - \lambda_n$.

As examples of such nonnegative measures $G_n$ $(n = 1, 2, \cdots)$, we may consider $G_n(\mathbf{x}) = 1$ $(\forall \mathbf{x} \in \mathcal{X}^n; \forall n = 1, 2, \cdots)$ with *countably infinite* source alphabet $\mathcal{X}$ (called the *counting measure* on $\mathcal{X}^n$) or the $n$-dimensional Lebesgue measure with *real* source alphabet $\mathcal{X}$. In particular, the case of the counting measure has the deep structural relationship with the fixed-length source coding problem, which will be elucidated in the next section.

**Remark 6.1** As will be easily seen from the proofs, even if we in turn replace the probability measure $P_{X^n}$ of the *null* hypothesis by nonnegative measures $F_n$ with $F_n(\emptyset) = 0$, both of Theorem 2.1 and Theorem 4.1 continue to hold with the due reinterpretation for probabilities as above. □

# 7 Hypothesis Testing and Fixed-Length Source Coding

Thus far, we have shown two key theorems (Theorem 2.1 and Theorem 4.1) concerning the general hypothesis testing. In this general setting, we can show also many other elegant systematic results on the hypothesis testing (as for the details, refer to Han [18]). In parallel with these systematic results, the corresponding many results in the general fixed-length source coding problem have been established (cf. Han [18, 20]). This correspondence is of very intrinsic nature not only at the technical level but also at the conceptual level, which can be made very transparent by introducing the *generalized* hypothesis testing problem as above. From this point of view, it turns out that all the theorems that hold in the fixed-length source coding problem can be regarded as forming a special class of those holding in the generalized hypothesis testing problem.

As an illustrative case, we will show that Theorem 2.1 of Han [20] immediately follows as a special case of Theorem 2.1 (in Section 2) with the counting measure $C_n(\mathbf{x}) \equiv 1$ ($\forall \mathbf{x} \in \mathcal{X}^n$) as the alternative hypothesis. To show this, let us first state the formal definition of the general fixed-length source coding problem. Let $\mathbf{X} = \{X^n\}_{n=1}^{\infty}$ be any general source with *countably infinite* source alphabet $\mathcal{X}$, and let $\mathcal{M}_n \equiv \{1, 2, \cdots, M_n\}$ be an integer set. Then, mappings $\varphi_n : \mathcal{X}^n \to \mathcal{M}_n$, $\psi_n : \mathcal{M}_n \to \mathcal{X}^n$ are called the *encoder* and the *decoder*, where we call $\varepsilon_n \equiv \Pr\{X^n \neq \psi_n(\varphi_n(X^n))\}$ the *error probability* of the fixed-length source coding. We denote the pair $(\varphi_n, \psi_n)$ with the error probability $\varepsilon_n$ by $(n, M_n, \varepsilon_n)$ (called a *code*). In the fixed-length source coding problem, we are interested in the prolem of determining the infimum $R_e(r|\mathbf{X})$ of achievable rates under asymptotic constraints of the form $\varepsilon_n \leq e^{-nr}$ ($r > 0$ is a prescribed constant) on the error probability $\varepsilon_n$. Formally, we define as follows.

**Definition 7.1** $R$ is called $r$-achievable if there exists a code $(n, M_n, \varepsilon_n)$ such that

$$\liminf_{n \to \infty} \frac{1}{n} \log \frac{1}{\varepsilon_n} \geq r,$$

$$\limsup_{n \to \infty} \frac{1}{n} \log M_n \leq R.$$

**Definition 7.2** *(The infimum of $r$-achievable fixed-length coding rates)*

$$R_e(r|\mathbf{X}) = \inf \{R \mid R \text{ is } r\text{-achievable}\}.$$

**Definition 7.3**

$$\sigma(R) = \liminf_{n \to \infty} \frac{1}{n} \log \frac{1}{\Pr\left\{\frac{1}{n} \log \frac{1}{P_{X^n}(X^n)} \geq R\right\}}. \tag{7.1}$$

With these definitions, the following general theorem has been established based on the entropy-spectrum argument which is a different version of the information-spectrum demonstrated in this paper.

**Theorem 7.1** (Han [18, 20]) Let $\mathbf{X} = \{X^n\}_{n=1}^{\infty}$ be a general source with *countably infinite* alphabet $\mathcal{X}$, then for any $r \geq 0$ we have

$$R_e(r|\mathbf{X}) = \sup_{R \geq 0} \{R - \sigma(R) \mid \sigma(R) < r\}, \tag{7.2}$$

where $R_e(0|\mathbf{X}) = 0$ ($r = 0$). □

Let us now show that Theorem 7.1 directly follows just by rewriting Theorem 2.1 with the counting measure $C_n(\mathbf{x}) \equiv 1 \quad (\forall \mathbf{x} \in \mathcal{X}^n)$ as the alternative hypothesis. Let this alternative hypothesis be denoted by $\mathbf{C} = \{C_n\}_{n=1}^{\infty}$. First, when we are given an acceptance region $\mathcal{A}_n \subset \mathcal{X}^n$ for a hypothesis testing, set $M_n = |\mathcal{A}_n|$ and we consider the encoder $\varphi_n : \mathcal{X}^n \to \mathcal{M}_n$ such that $\varphi_n$ maps in the one-to-one manner all the elements of $\mathcal{A}_n$ into $\mathcal{M}_n$ in the order of $1, 2, \cdots$, and maps all the elements of $\mathcal{A}_n^c$ into $1 \in \mathcal{M}_n$, where the decoder $\psi_n : \mathcal{M}_n \to \mathcal{X}^n$ is the inverse mapping of $\varphi_n|_{\mathcal{A}_n}$. Then, it is obvious that
$$\mathcal{A}_n = \{\mathbf{x} \in \mathcal{X}^n \mid \psi_n(\varphi_n(\mathbf{x})) = \mathbf{x}\},$$
which means that the first kind of error probability $\mu_n = \Pr\{X^n \notin \mathcal{A}_n\}$ for the hypothesis testing coincides with the error probability $\varepsilon_n$ for the fixed-length source coding. We notice that this kind of correspondence between hypothesis testings and fixed-length source codings becomes the one-to-one mapping if we indifferently identify all the codes $(n, M_n, \varepsilon_n)$ which have the same set $\mathcal{A}_n = \{\mathbf{x} \in \mathcal{X}^n \mid \psi_n(\varphi_n(\mathbf{x})) = \mathbf{x}\}$ of the elements of $\mathbf{x} \in \mathcal{X}^n$ that can be correctly decoded under fixed-length source coding. On the other hand, the second kind of error probability $\lambda_n$ under the counting measure $C_n$ can be written as
$$\begin{aligned}\lambda_n &= C_n(\mathcal{A}_n) = |\mathcal{A}_n| = M_n \\ &= e^{nr_n},\end{aligned} \qquad (7.3)$$
where
$$r_n = \frac{1}{n}\log M_n.$$
Then, under this correspondence it follows from (7.3) that
$$\liminf_{n \to \infty} \frac{1}{n} \log \frac{1}{\lambda_n} = -\limsup_{n \to \infty} r_n,$$
which means that $R$ is $r$-achievable for (generalized) hypothesis testing if and only if $-R$ is $r$-achievable for fixed-length source coding. Thus, from Definition 7.1 Definition 7.2 and Definition 2.1 Definition 2.2, we have the following equation connecting $B_e(r|\mathbf{X}||\mathbf{C})$ to $R_e(r|\mathbf{X})$:
$$B_e(r|\mathbf{X}||\mathbf{C}) = -R_e(r|\mathbf{X}) \quad (\forall r > 0). \qquad (7.4)$$
Next, since we are considering the counting measure $C_n$ as the alternative hypothesis, the probability appearing on the righ-hand side of (2.1) defining $\eta(R)$ is written as
$$\begin{aligned}&\Pr\left\{\frac{1}{n}\log \frac{P_{X^n}(X^n)}{P_{\overline{X}^n}(X^n)} \leq R\right\} \\ &= \Pr\left\{\frac{1}{n}\log \frac{P_{X^n}(X^n)}{C_n(X^n)} \leq R\right\} \\ &= \Pr\left\{\frac{1}{n}\log P_{X^n}(X^n) \leq R\right\} \\ &= \Pr\left\{\frac{1}{n}\log \frac{1}{P_{X^n}(X^n)} \geq -R\right\}.\end{aligned}$$

Therefore, we have $\eta(R) = \sigma(-R)$ by the definition (7.1) of $\sigma(R)$, i.e.,

$$\sigma(R) = \eta(-R). \tag{7.5}$$

Then, Theorem 2.1 with the counting measure as the alternative hypothesis together with (7.4) yields

$$\begin{aligned} R_e(r|\mathbf{X}) &= -B_e(r|\mathbf{X}||\mathbf{C}) \\ &= -\inf_R \{R + \eta(R) \mid \eta(R) < r\} \\ &= \sup_R \{-R - \eta(R) \mid \eta(R) < r\}. \end{aligned}$$

As a consequence, if we replace $R$ by $-R$ and use (7.5), it is concluded that

$$R_e(r|\mathbf{X}) = \sup_{R \geq 0} \{R - \sigma(R) \mid \sigma(R) < r\}.$$

This is nothing but Theorem 7.1 on the fixed-length source coding.